\documentclass[sigconf, nonacm]{acmart}
\usepackage{booktabs}
\usepackage{array, makecell}
\usepackage{subcaption}
\usepackage{mathrsfs}
\usepackage{multirow}
\usepackage{soul,color}
\usepackage{algorithm}
\usepackage{algpseudocode}
\usepackage{tabularx}
\usepackage{dblfloatfix}
\usepackage{wrapfig}
\AtBeginDocument{%
  \providecommand\BibTeX{{%
    \normalfont B\kern-0.5em{\scshape i\kern-0.25em b}\kern-0.8em\TeX}}}

\acmSubmissionID{994}

\begin{document}


\title[Data-Driven Optimization for Police Districting in South Fulton, Georgia]{Data-Driven Optimization for Police Districting \\in South Fulton, Georgia}

\author{Shixiang Zhu}
\email{shixiang.zhu@gatech.edu}
\affiliation{%
  \institution{Georgia Institute of Technology}
  \streetaddress{North Ave NW}
  \city{Atlanta}
  \state{Georgia}
  \postcode{30332}
}

\author{Alexander W. Bukharin}
\email{abukharin3@gatech.edu}
\affiliation{%
  \institution{Georgia Institute of Technology}
  \streetaddress{North Ave NW}
  \city{Atlanta}
  \state{Georgia}
  \postcode{30332}
}

\author{Le Lu}
\email{helenlule@gatech.edu}
\affiliation{%
  \institution{Georgia Institute of Technology}
  \streetaddress{North Ave NW}
  \city{Atlanta}
  \state{Georgia}
  \postcode{30332}
}

\author{He Wang}
\email{he.wang@isye.gatech.edu}
\affiliation{%
  \institution{Georgia Institute of Technology}
  \streetaddress{North Ave NW}
  \city{Atlanta}
  \state{Georgia}
  \postcode{30332}
}

\author{Yao Xie}
\email{yao.xie@isye.gatech.edu}
\affiliation{%
  \institution{Georgia Institute of Technology}
  \streetaddress{North Ave NW}
  \city{Atlanta}
  \state{Georgia}
  \postcode{30332}
}

\renewcommand{\shortauthors}{Zhu, et al.}

\begin{abstract}
We redesign the police patrol beat in South Fulton, Georgia, in collaboration with the South Fulton Police Department (SFPD), using a predictive data-driven optimization approach. Due to rapid urban development and population growth, the existing police beat design done in the 1970s was far from efficient, which leads to low policing efficiency and long 911 call response time. We balance the police workload among different city regions, improve operational efficiency, and reduce 911 call response time by redesigning beat boundaries for the SFPD. We discretize the city into small geographical atoms, which correspond to our decision variables; the decision is to map the atoms into ``beats'', the basic unit of the police operation. We first analyze workload and trend in each atom using the rich dataset, including police incidents reports and U.S. census data; We then predict future police workload for each atom using spatial statistical regression models; Lastly, we formulate the optimal beat design as a mixed-integer programming (MIP) program with continuity and compactness constraints on the beats' shape. The optimization problem is solved using simulated annealing due to its large-scale and non-convex nature. The simulation results suggest that our proposed beat design can reduce workload variance among beats significantly by over 90\%. 
\end{abstract}

%
%




\maketitle

\section{Introduction}
\label{sec:introduction}

The City of South Fulton, Georgia, was recently established in May 2017 from previously unincorporated land outside Atlanta. It is now the third-largest city in Fulton County, Georgia, and serves a population of over 98,000, among which 91.4\% are black, or African American  \cite{UScensusbureau}. South Fulton is a historic area renowned for its art and activism. Despite this, the city has often faced the challenge of climbing crime rates and long police response times. In a 2019 survey, 46.48\% of residents responded that they do not feel safe in South Fulton. In the same year, the South Fulton City Council made it clear that their number one priority was to make South Fulton safer \cite{SFStrategicPlan}. 

\begin{figure}[!t]
\centering
\includegraphics[width=.7\linewidth]{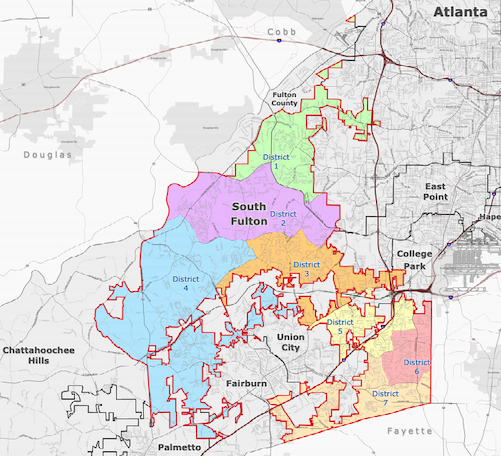}
\caption{City-wide police district map of South Fulton, GA. There were 7 beats, which was initially designed in the 1970s. The city boundary is highly irregular which requires intricate design of police beats.}
\label{fig:city-map}
\vspace{-0.2in}
\end{figure}

The South Fulton Police Department (SFPD) is the main policing force in the city. From 2019 to early 2020, our team worked with the SFPD to improve their police operation efficiency. Our project specifically focused on redistricting beat configurations (by completely re-drawing the beat boundaries and changing the number of beats), aiming to rebalance SFPD officers' \emph{workload} (total amount of working time). The initial analysis identified that workload unbalance among different areas of the city was caused by an outdated beat design that had not been changed for over five decades; the inefficient beat design, in turn, lead to long 911 call response time in some areas.

Previously, the police operation of South Fulton was according to seven
police \emph{beats}, which divide the city geographically as shown in Figure~\ref{fig:city-map}. 117 police personnel were allocated to the beats for patrolling and responding to the 911 calls \cite{SouthFultonWebsite}. Typically, at each shift, one response unit (usually a police car with one to two officers) answers all the 911 calls that occurred in a certain beat. If the response unit is busy handling another incident, nearby available response units may be dispatched by the operator to answer the call. 

\begin{figure}
\centering
\includegraphics[width=.8\linewidth]{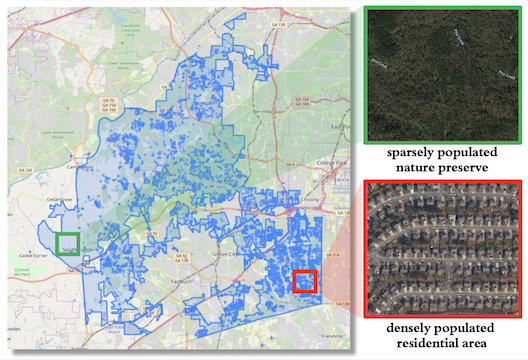}
\vspace{-0.1in}
\caption{Distribution of 911 calls-for-services requests in South Fulton, GA. Blue shaded area is the city limit of South Fulton. Blue dots are locations of requests. The requests are unevenly distributed among different regions.}
\label{fig:911call-distribution}
\vspace{-0.2in}
\end{figure}

The most recent South Fulton police beat redesign occurred in the 1970s -- almost five decades ago. Since then, the area (which eventually became the City of South Fulton) has undergone tremendous urban growth that drastically changes its landscape. The U.S.\ Census Bureau estimated that South Fulton's population has increased by 13.7\% from 2010 to 2018 \cite{UScensusbureau}. 
The city's rapid development has led to a significant increase in police workload, exacerbated by the difficulty in officer recruitment and retention faced by the SFPD. 
Moreover, demographic and traffic pattern changes also create an unbalanced workload among different regions. Figure~\ref{fig:911call-distribution} shows the distribution of 911 calls, recorded by real 911-call reports provided by SFPD from 2018 to 2019. The figure shows some beats faced a significantly higher workload than others. For example, police officers in the city's southeastern area respond to more calls than those in the western region.

Since the seminal work by R.\ Larson and others \cite{larson1972urban, larson1981urban}, researchers have recognized that beat configuration may significantly impact police response time to 911 calls and operational efficiency. In particular, the area and shape of beats determine the workload and travel time in that beat. Hence, it is critical to design the boundaries of beats to balance the workload.

\begin{figure}[!t]
\centering
\includegraphics[width=1.\linewidth]{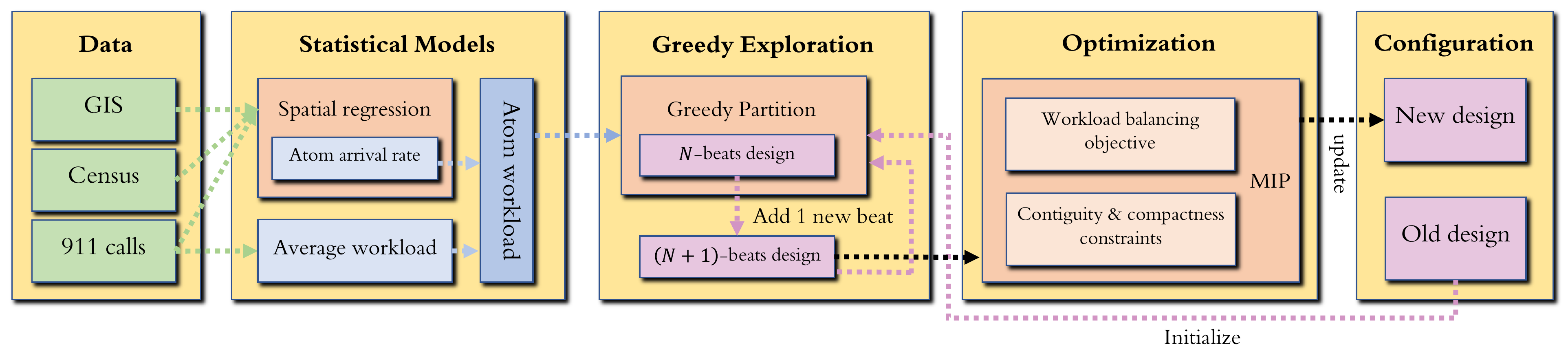}
\vspace{-0.3in}
\caption{An illustration for data-driven optimization framework of police beat redesign.}
\label{fig:workflow}
\vspace{-0.2in}
\end{figure}

\vspace{.05in}
\noindent{\bf Outline.}
We redesigned the police patrol beats in the City of South Fulton using a data-driven optimization approach. The outline of our approach is summarized in Figure~\ref{fig:workflow}.
Our objective is to balance police workload in each beat by redrawing beat boundaries. First, we divided the geographical areas of the city into a large number of ``atoms''. Then, we estimated the workload in each atom using police reports data and census data, including population and socio-economic factors. These steps are described in Sec.~\ref{sec:data_sources} and \ref{sec:data-preprocessing}.
Based on the workload estimation, we developed statistical models to predict police workload in the next few years (Sec.~\ref{sec:workload-prediction}).  
We then formulate the beat redesign problem as a clustering problem: each beat is formed with a cluster of atoms. This clustering problem is formulated and solved using mixed-integer programming (MIP), where
the objective function is a metric of workload unbalance (defined as the workload variance across all beats). We also impose constraints that require beats to be contiguous and compact so that they are not irregularly shaped. The problem formulation is described in Sec.~\ref{sec:beat_redesign_opt}. To tackle the computational complexity of solving a large-scale optimization problem, we developed a simulated annealing based approach with efficient solution exploration. We also study the districting with the different number of beats and find the optimal number of beats with the highest cost-effectiveness. Numerical results (Sec.~\ref{sec:results}) show that our proposed beat design can reduce workload variance among different regions by over 90\%. In January 2020, together with the SFPD, we presented our final redesign plan to the South Fulton city council, which was officially approved for implementation. 

\vspace{.05in}
\noindent{\bf Contribution.}
Our work proposes a new data-driven framework that integrates data, statistical prediction, and optimization in the context of police beat design. Previous works in the predictive policing literature tend to focus on only the prediction aspect. The operations research literature often studies police zone design based on analyzing stochastic models without explicitly considering data sources. We take advantage of the availability of abundant data and adopt a new data-driven approach: the workload and other important parameters for optimization are estimated and predicted from data. From a methodological perspective, 
we use geo-spatial atoms to define city boundaries and police beat boundaries. This approach enables accurate workload prediction by correlating historical police data with the census data and beat design optimization. 
 
Our project also had a significant societal impact and directly improved the police operations of the SFPD and the safety of residents in South Fulton. It is worth mentioning that although we focus on the study of police beat redesign in South Fulton, our method can be applied to other cities facing similar issues.

\vspace{.05in}
\noindent{\bf Related work.} 
Police districting (designing beats or zones) is a classical problem studied in operations research dating back to the 1970s (see the seminal work \cite{larson1972urban} and the surveys by \cite{chaiken1972methods,green2004anniversary} for reviews). \cite{gass1968division} is one of the earliest works that study optimal beat allocation using integer programming. \cite{bammi1975allocation} considers the beat allocation problem to minimize response time for police service calls. In particular, the paper also considers overlapping beats, where multiple patrol officers share one patrol area. \cite{chelst1979estimating,Larson1974} use queueing models to estimate travel time. In particular, our proposed data-driven model includes the travel time in the workload calculation.
\cite{larson1974illustrative} introduces a heuristic approach to the design of beats with implementation in Boston. \cite{bodily1978police} considers fairness issues of police zone design. We remark that most classical works rely on analyzing stochastic models for police workload estimation, which usually requires stringent assumptions, e.g., calls arrive according to homogeneous Poisson processes (with the notable exception of \cite{Larson1974}). Here, rather than obtained from stochastic models, we take advantage of the availability of abundant data and adopt a data-driven approach: the workload and other important parameters for optimization are estimated and predicted from data.

There is also a large body of works on other types of geographical districting problems, such as political districting. This includes the pioneering work \cite{GaNe1970} studies political districting using integer programming. Their method is extended by \cite{Sh2009} for other geographical districting problems. A few other works \citep{WeHe2963, Mi1967, Mo1973, Mo1976, Vi1961, Damico2002} apply meta-heuristics (e.g., genetic algorithms, simulated annealing) to geographical districting, which usually lack optimality guarantees. Geographic districting often include criteria such as contiguity \citep{Gr1985, Mi1967, GaNe1970, Na1972, Me1972, Vi1961} and compactness \citep{GaNe1970, Ni1990, Yo1988}, which are also important in the police zone design context. However, political districting has different considerations than police districting.

In the last decade, we have seen the rise of predictive policing, i.e., the use of mathematical and statistical methods in law enforcement to predict future criminal activity based on past data. Its importance has been even recognized by Time magazine that in November 2011 named predictive policing as one of the 50 best inventions of 2011 \citep{Grossman2011}.
The RAND Corporation and the National Institute of Justice of the United States (NIJ) also acknowledge the need for taking a step forward and developing explicit methodologies and tools to take advantage of the information provided by predictive policing models to support decision makers in law enforcement agencies \citep{rios2020optimal}.
We remark that most classical works do not leverage the historical operational data and rely on analyzing stochastic models for police workload estimation, which usually requires stringent assumptions, e.g., calls arrive according to homogeneous Poisson processes. Here, rather than obtained from stochastic models, we take advantage of the availability of abundant data and adopt a data-driven approach: the future workload and other essential parameters for optimization are estimated and predicted from data.


\vspace{-.05in}
\section{Data}
\label{sec:data_sources}

We start by describing the various sources of data used for South Fulton police beats reconfiguration, including 911 calls-for-service reports, geographical data of the city, and the socio-economic data collected by the American Community Survey (ACS) from the U.S.\ Census Bureau.

\vspace{.05in}
\noindent{\bf 911 calls-for-service data.}
The SFPD provides comprehensive 911-call reports between May 2018 to April 2019, which contains 69,170 calls in total (Figure~\ref{fig:911call-distribution}). The recorded 911 calls cover more than 600 categories of incidents, including assaults, terrorist threats, domestic violence, robbery, burglary, larcenies, auto-thefts, etc. 
These reports are generated by mobile patrol units in the city, which handle 911 calls 24/7. 
Teams of \emph{response units} (police cars and officers) are assigned to patrol city streets, and answer calls for service. 
When a 911 call for a traffic incident comes in at the \emph{call time}, a new incident record will be created at the dispatch center, and the call location will be recorded. The operator assigns an officer to handle the call. The unit arrives at the scene and starts the investigation. Once the police complete the investigation and clear the incident, the police report will be closed and record the \emph{clear time}. The time interval that it takes police to process the call between the call time and the clear time is called \emph{processing time}.
The police workload is calculated using both the geolocation data and 911 call processing time data (The calculation method which will be discussed in more detail in Sec.~\ref{sec:data-preprocessing}). The geolocation consists of the GPS location of reported incidents. From the geographical data of South Fulton, we are also able to identify which beat each incident is located. 

\vspace{0.05in}
\noindent{\bf GIS data \& beat configuration.}
Geographic information system (GIS) data contain the geographical information of the city's and beats' boundaries, which are extracted from Fulton County Special Services District digest parcel data \cite{SFGIS}.
Geographically, the city boundary of South Fulton is quite irregular with jagged edges, holes, and disconnected segments (Figure~\ref{fig:city-map}). This irregularity is due to the formation of the City of South Fulton, with the city being a new combination of all the unincorporated land in southwest Fulton County.
Currently, there are seven beats in the City of South Fulton.
As shown in figure ~\ref{fig:city-map}, beat (district) 1, 2, 3, and 4 include larger areas that are relatively compact, while remaining beats contain smaller scattered areas. The irregular shape of the city brings difficulty to police officers while reaching locations of requests and patrolling. 
Moreover, as the busiest airport globally, the Hartsfield-Jackson Atlanta International Airport is situated east of the city, which significantly adds to the city's workload disparity.

\begin{figure}[!t]
\centering
\begin{subfigure}[h]{0.32\linewidth}
\includegraphics[width=\linewidth]{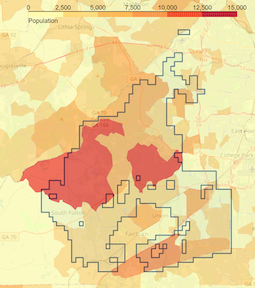}
\caption{Population}
\label{fig:census-population}
\end{subfigure}
\begin{subfigure}[h]{0.32\linewidth}
\includegraphics[width=\linewidth]{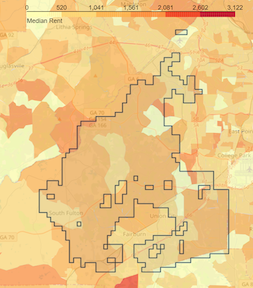}
\caption{Median rent}
\label{fig:census-rent}
\end{subfigure}
\begin{subfigure}[h]{0.32\linewidth}
\includegraphics[width=\linewidth]{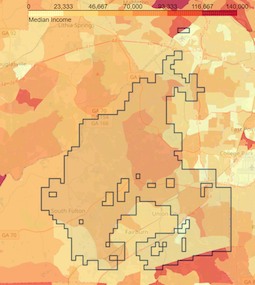}
\caption{Median income}
\label{fig:census-income}
\end{subfigure}
\begin{subfigure}[h]{0.32\linewidth}
\includegraphics[width=\linewidth]{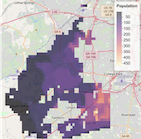}
\caption{Population}
\label{fig:atom-population}
\end{subfigure}
\begin{subfigure}[h]{0.32\linewidth}
\includegraphics[width=\linewidth]{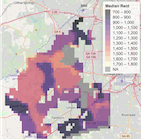}
\caption{Median rent}
\label{fig:atom-rent}
\end{subfigure}
\begin{subfigure}[h]{0.32\linewidth}
\includegraphics[width=\linewidth]{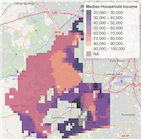}
\caption{Median income}
\label{fig:atom-income}
\end{subfigure}
\vspace{-0.1in}
\caption{(a-c): Raw data for demographic factors of South Fulton, GA in 2019, from American Community Survey, organized by census blocks. (d-f): Corresponding atomized census data of South Fulton, GA, in 2019.}
\label{fig:tau}
\vspace{-0.2in}
\end{figure}

\vspace{0.05in}
\noindent{\bf Census data.}
The American Community Survey (ACS) collected by the U.S.\ Census Bureau provides comprehensive information about the population, demographic, and economic status of different Georgia areas. Unlike the census, which takes place every ten years, the ACS is conducted once per year. 
Some demographic factors are useful in predicting future workload (by correlating the city's socio-economic profile with the workload). These factors contain essential information about the development and economic growth of the city. 
Besides, census data is organized by census blocks, as shown in Figure~\ref{fig:tau} (a - c), which is also different from the geographical atoms we consider for our study.
In our discussions with the SFPD, we selected eight most influential factors that play a vital role in determining the police workload, such as population and median rent, school enrollment, and the average year structures were built. The full list of census factors we are considering has been shown in Table~\ref{tab:regr-table}.

\vspace{-.05in}
\section{Data preprocessing}
\label{sec:data-preprocessing}

In this section, we describe three key steps in data preprocessing before performing the beat design. In particular, we need to address the following challenge in using the data: how to align time resolution and spatial resolution from the raw data with what we need in the design.

\begin{figure}
\centering
\includegraphics[width=.8\linewidth]{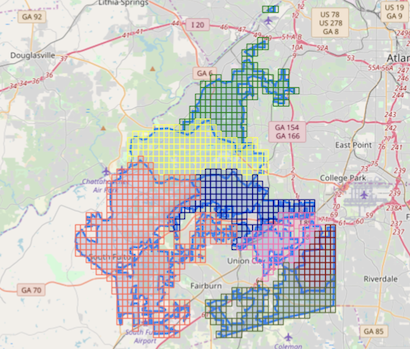}
\caption{South Fulton region is partitioned into 1,187 square geographical atoms. Color indicates the beat membership.}
\label{fig:city-atom}
\vspace{-0.1in}
\end{figure}

\vspace{0.05in}
\noindent{\bf Geographical atoms.}
To accurately capture changing demographics and determine the new boundaries for each police beat, we define high-resolution geographical atoms 
by creating artificial polygons of identical size as our geographical atoms to break up the city. The optimal beat design can be found by aggregating multiple adjacent polygons. 
The size of geographical atoms is essential to our design's performance since it determines the number of variables in the optimization and the precision of the workload estimation. There is a trade-off between computational efficiency and model accuracy in determining the size of geographical atoms. If the atoms' size is too large, then we are unable to capture community demographics accurately; if the size of the atoms is too small, then the problem will become computationally intractable. 
After rounds of discussion with the SFPD, we decide atoms to be a square area with a side length of 0.345 miles, roughly the city block size.
This allows us to estimate the local workload accurately while resulting in a reasonable number of decision variables in our optimization problem. 
The atomized map of the city was generated by intersecting the city boundary with a grid of atoms, resulting in a new grid of 1,187 geographical atoms, as shown in Figure~\ref{fig:city-atom}. The police workload estimation and prediction will be performed based on these predefined geographical atoms.
Formally, let $i\in \mathscr{I}=\{1,\dots,I\}$ denote the $i$-th atom and $k \in \mathscr{K} = \{1,\dots, K\}$ denote the $k$-th beat in our design. Let the binary decision variable $d_{ik} \in \{0, 1\}$ denote whether or not atom $i$ is assigned to beat $k$. 
A particular beat design is a unique graph partition determined by a matrix $D=\{d_{ik}\} \in \{0,1\}^{I \times K}$. For each $i$, it satisfies $\sum_{k = 1}^{K} d_{ik} = 1$. 
Given the beat design $D$, the set of atoms assigned to beat $k$ is denoted by $\mathscr{I}_k(D) = \{i: d_{ik} = 1\} \subseteq \mathscr{I}$. %
%
Figure~\ref{fig:city-atom} also shows the discretization of the existing beat configuration, where atoms with the same color represent a police beat. 


\vspace{0.05in}
\noindent{\bf Census data atomization.}
A major challenge for estimating the socio-economic data for each geographical atom using census data is the inconsistency between census blocks and geographical atoms, where, as shown in Figure~\ref{fig:tau} (a-c), census blocks usually have a much larger area than geographical atoms. Here we need to perform a spatial interpolation to align the census data with our geographical atoms. 
Specifically, we assume the census data, such as population, in the same geographical atom, are evenly distributed. The data of each census factor in a geographical atom can be estimated by proportionally dividing the value in the census block where the atom falls into. The weight of the portion that an atom takes from a census block can be measured by the proportion between their areas. As shown in Figure~\ref{fig:tau} (d - f), the census data collected by census blocks have been discretized into geographical atoms. 
Given historical census data in the month $\ell \in [L - L_0, L]$, where $L$ and $L_0$ denote the last month and the time span of the historical data, respectively.
%
The preprocessed census data is denoted as a tensor $X = \{x_{i \ell m}\} \in \mathbb{R}^{I \times L_0 \times M}$,
where each entry $x_{i \ell m}$ indicates the value of the census factor $m \in \mathscr{M} = \{1, \dots, M\}$ in atom $i$ and month $\ell$.


\vspace{0.05in}
\noindent{\bf 911 calls-for-service data preprocessing.}
We estimate the police workload for each geographical atom using the 911 calls-for-service dataset. 
The workload of each 911 call is evaluated by its processing time, i.e., the total time that the police spend on traveling and the investigation.
We calculate the workload by two steps: (1) count the number of 911 calls occurred in the $i$-th atom in $\ell$-th month, denoted as $N_{i\ell}$; (2) estimate the total workload for the $i$-th atom in the $\ell$-th month by multiplying $N_{i\ell}$ by the average processing time, denoted as $w_{i\ell}$. The count of 911 calls will be further used as the predictor in our spatial regression model, which will be discussed in Sec.~\ref{sec:workload-prediction}.




\vspace{-.05in}
\section{Police workload Prediction}
\label{sec:workload-prediction}

Predicting the police workload is particularly challenging.
Although we assumes the call arrival rates are time-homogeneous,
we observe in the actual data that the call arrival rates have a significant seasonality pattern and yearly trend, as well as correlation over adjacent geographical areas.
Therefore, we propose a spatio-temporal model to predict future call arrival rates.
We assume that the call arrival rate $\lambda_{i\ell}$ for atom $i$ in a month $\ell$ is a constant. Thus, each beat is a homogeneous Poisson process with rate $\lambda_{i\ell}$. The arrival rates $\Lambda = \{\lambda_{i\ell}\} \in \mathbb{R}_+^{I \times L_0}$ can be approximated by $N_{i\ell}$, where $L_0 = 12$. 
We learned from the SFPD that, the occurrence of 911 calls is highly correlated with population and economic status of the beat and its neighborhood. 
We predict the arrival rate $\lambda_{i\ell}$ in the future month $\ell = L + t, t = 1, 2, \dots$ using a linear model that regresses the arrival rate to other endogenous variables (arrival rates in other beats) and exogenous factors (demographic factors). 
As shown in Table~\ref{tab:regr-table}, we consider $M=8$ demographic factors, which are statistically verified to be good predictors, including population, education level, and household income. Specifically, we use the spatially lagged endogenous regressors \citep{Rosen1974} defined as
\begin{equation*}
  \lambda_{i\ell} = \sum_{(i,j)\in\mathscr{A}} \alpha_{ij} \lambda_{j\ell} + \beta_0 \lambda_{i,\ell-1} + \sum_{t=1}^{p} \boldsymbol{\beta}_t^\intercal X_{i,\ell-t} + \epsilon_{i}, \quad \forall \ell \in [L-L_0, L],
  \label{eq:lr-lam}
\end{equation*}
where $p$ is the total number of past months of data that we consider for fitting the regressor, which in our case was 1. The adjacency matrix $A = \{\alpha_{ij}\} \in \mathbb{R}^{I \times I}$ specifies adjacency relationships between atoms. The temporal coefficient $\beta_0 \in \mathbb{R}$ specifies the influence of the last month. The coefficient $\boldsymbol{\beta}_t \in \mathbb{R}^{M}, \forall 1 \le t \le p$ specifies correlations with census factors and  error term $\epsilon_{i}$ are spatially correlated. The set of adjacency pairs is defined by $\mathscr{A} = \{(i,j): i,j\text{ are adjacent in }\mathcal{G};\ i, j \in \mathscr{I}\}$. The graph $\mathcal{G}$ is given by associating a node with every atom and connecting two nodes by an edge whenever the corresponding atoms are geographically adjacent. Here, we capture the spatial correlation between data using the standard spatial statistics approach, by assuming $\epsilon_i$ to be spatially correlated with correlation depending distance between two locations \cite{ripley2005spatial}. 


\begin{table}[t]
\caption{Variables used for workload prediction}
\vspace{-0.15in}
\label{tab:regr-table}
\begin{center}
\begin{small}
\begin{sc}
\resizebox{0.48\textwidth}{!}{%
\begin{tabular}{ccc}
\toprule[1pt]\midrule[0.3pt]
     Predictor & Regression Coefficient & p-value \\ [0.5ex] 
     \hline
    Population & 439.558 & 0.007 \\ 
    Number of housing units & 158.440 & 0.019 \\
     School Enrollment & 79.236 & 0.008 \\
    Median Household Income & 59.420 & 0.000 \\
     Median Number of Rooms & -10.560 & 0.006 \\
    Median Age & -7.421 & 0.001 \\
    Median House Price & -16620 & 0.000 \\
    Average Year Built & 170.140 & 0.003  \\
\midrule[0.3pt]\bottomrule[1pt]
\end{tabular}
}
\end{sc}
\end{small}
\end{center}
\end{table}

\section{Beat Redesign Optimization}
\label{sec:beat_redesign_opt}

In this section, we introduce our objective and solution methods to the beat redesign optimization problem. We develop an optimization framework to shift beat boundaries, where artificial geographical atoms were assigned to beats while balancing the workload. We formulate this problem as minimizing the workload variance by reconfiguring the beat plan with constraints, including the continuity and compactness of beats. 

\subsection{Objective}
\label{sec:objective}

Our goal is to shift beat boundaries and make inter-beat workload distribution even. 
Based on the discussion with the police, we choose the objective function as the workload variance among different zones, 
which quantitatively measure the police workload imbalance between zones from a macro view.
The objective of this problem can be formulated as minimizing the inter-beat workload variance $Z(D)$ given a beat design $D$:
\begin{equation}
\begin{split}
\underset{D}{\text{minimize}} &\quad  Z(D) := \sum_{k=1}^{K} \left( w_{k\ell}(D) - \frac{\sum_{\kappa=1}^K w_{\kappa\ell}(D)}{K} \right)^2\\
\mbox{subject to} 
&\quad \sum_{k = 1}^{K} d_{ik} = 1,\quad \forall i\\
&\quad \mbox{contiguity and compactness for each beat}.
\end{split}
\label{eq:opt-objective-quadratic-1}
\end{equation}
Recall that the matrix $D = \{d_{ik}\} \in \{0,1\}^{I \times K}$ represents decision variables, where binary variable $d_{ik} \in \{0,1\}$ indicates whether or not geographical atom $i$ is assigned to beat $k$; and $w_{k\ell}(D) = \sum_{i \in \mathscr{I}_k(D)} w_{i\ell} $ represents the total workload in beat $k=1,\cdots,K$ in month $\ell$.
The variance is a quadratic function of the workload in each beat, which implies that the objective function is convex with respect to the decision variables. A smaller variance indicates a more balanced inter-beat police workload.
The constraints will be explicitly defined in Sec.~\ref{sec:constraints}.

\subsection{Compactness and contiguity constraints}
\label{sec:constraints}

In addition to balancing the police workload, it is desirable that the beat shapes are contiguous and compact. 
In fact, the police never used a quantitative measure of compactness to declare the plans unsuitable. Instead, the police have simply disallowed plans with long and thin or snakelike districts. In other words it appears that the police have evaluated compactness only visually. 
Since it is not obvious how to determine an acceptable compact design, we choose to minimize the workload variance based on the discussion with the police; 
but it should be understood that compactness is in reality a loose constraint rather than an objective.
Therefore, we formulate the \emph{contiguity} and \emph{compactness} criteria as a set of linear constraints \cite{GaNe1970, Ni1990, Shirabe2009, Yo1988} by introducing additional variables: $f_{ijk}$ is the flow from atom $i$ to atom $j$ in beat $k$; $h_{ik}$ equals to 1 if atom $i\in\mathscr{I}$ is selected as a sink in beat $k\in\mathscr{K}$, otherwise 0; $q$ is the maximum beat capacity.
Hence, there are 21,170,145 variables with 63,421,410 constraints in total.

\vspace{.05in}
\noindent\textbf{Contiguity constraints.} Contiguity constraints are imposed on each beat using the flow method \cite{Shirabe2009}. For each beat $k$, there is a flow $f_{ijk}$ on the graph, where $f_{ijk}$ denotes flow from $i$ to $j$. Each beat has a hub vertex whose net flow is at most the number of vertices in the beat, less one. Each other vertex in the beat has a net flow of at most $-1$. This ensures that there is a path of positive flow from any vertex in the beat to the hub, implying contiguity.

Specifically, constraints \eqref{eq:opt-con2} represent the net outflow from each beat. The two terms on the left indicate, respectively, the total outflow and total inflow of atom $i$. If atom $i$ is included in beat $k$ but is not a sink, then we have $d_{ik}=1$, $h_{ik}=0$, and thus atom $i$ must have supply $\ge1$. If atom $i$ is included in beat $k$ and is a sink, then we have $d_{ik}=1$, $h_{ik}=1$, and thus atom $i$ can have demand (negative net outflow) $\le q-1$. If atom $i$ is not included in beat $k$ and is not a sink, then we have $d_{ik}=0$, $h_{ik}=0$, and thus atom $i$ must have supply 0. If atom $i$ is not included in beat $k$ but is a sink, then we have $d_{ik}=0$, $h_{ik}=1$, and the rest of $d_{\cdot k}$ are forced to be 0, that is, no atoms are selected. Constraints \eqref{eq:opt-con3} specify the number of atoms that can be used as sinks. Constraints \eqref{eq:opt-con4} ensure that each beat must have only one sink. Constraints \eqref{eq:opt-con5} ensure that there is no flow into any atom $i$ from outside of beat $k$ (where $d_{ik}=0$), and that the total inflow of any atom in beat $k$ (where $d_{ik}=1$) does not exceed $q-1$. Constraints \eqref{eq:opt-con6} make sure unless a atom $i$ is included in beat $k$, the atom $i$ cannot be a sink in beat $k$. Constraints \eqref{eq:opt-con7} and \eqref{eq:opt-con8} ensure that there are no flows (inflows and outflows) between different beats which forces eligible contiguity.
\begingroup
\allowdisplaybreaks
\begin{subequations}
\begin{align}
  \sum_{(i, j) \in \mathscr{A}} f_{ijk} - \sum_{(i, j) \in \mathscr{A}} f_{jik} & \ge d_{ik} - q h_{ik}, & \forall i, k, \label{eq:opt-con2}\\
  \sum_{k}^{K} \sum_{i}^{N} h_{ik} & = K, & \label{eq:opt-con3}\\
  \sum_{i}^{N} h_{ik} & = 1, & \forall k, \label{eq:opt-con4}\\
  \sum_{(i, j) \in \mathscr{A}} f_{jik} & \le (q-1) d_{ik}, & \forall k, \label{eq:opt-con5}\\
  h_{ik} - d_{ik} & \le 0, & \forall i,k, \label{eq:opt-con6}\\
  f_{ijk} + f_{jik} & \le (q-1) d_{ik}, & \forall i,k, \label{eq:opt-con7}\\
  f_{ijk} + f_{jik} & \le (q-1) d_{jk}, & \forall j,k, \label{eq:opt-con8}\\
  d_{ik}, h_{ik} & \in \{0, 1\}, & \forall i, k, \label{eq:opt-con9}\\
  f_{ijk} & \ge 0, & \forall i, j, k, \label{eq:opt-con10}
\end{align}
\label{eq:opt-con-contiguity}
\end{subequations}
\endgroup

\noindent\textbf{Compactness constraints.} Compactness is defined as geographical compactness with distance compactness and shape compactness \cite{Ni1990, Yo1988}. For distance compactness, a district is feasible only if the distance between population units must be less than a specified upper bound. For shape compactness, a district is feasible only if the square of the distance's maximum diameter divided by the district's area must be less than another upper bound \cite{GaNe1970}. 
  
Following the existing literatures, we add two additional linear constraints \eqref{eq:opt-con11}, \eqref{eq:opt-con12} to ensure the compactness of beats. For each atom $i$, let $A_i$ be the area of $i$, and for each pair of atoms $i$ and $j$, let $l_{ij}$ be the square of the distance between the centroids of the beats. We also have a parameter $c_1, c_2>0$ controlling the degree of compactness.
\begin{subequations}
\begin{align}
  l_{ij}e_{ijk} & \leq c_1, & \forall i, j, k, \label{eq:opt-con11}\\
  l_{ij}e_{ijk} & \leq c_2 \sum_{i=1}^{K} d_{ik} A_i,  & \forall i, j, k, \label{eq:opt-con12}
\end{align}
\label{eq:opt-con-compact}
\end{subequations}

\subsection{Heuristic approximation}
\label{sec:methods}

Three methods were discussed in our experiments to search for optimal police beat design. The greedy algorithm serves to generate new beats iteratively and confirms the optimal number of beats for the future redesign. Following the greedy redesign, we adopt a heuristic optimization approach to find the beat design in contrast to the mixed-integer programming (MIP) approach.


\vspace{.05in}
\noindent{\bf Greedy search.}
To determine the optimal number of beats in the final design, we perform an iterative greedy algorithm, which attempts to generate new beat greedily for the design for each iteration while preserving the original structure of the existing beat as much as possible. Intuitively, more beats may result in a more balanced workload distribution. However, the manpower of the SFPD and resources of the South Fulton City Council are limited. It is unrealistic to deploy such a design with a large number of beats. 
Hence, we adopt the Greedy algorithm to explore the optimal number of beats in our design.  
The procedure for ``Greedily'' creating new beat designs is demonstrated as follows.

For the $n$-th iteration, we define $D_n$ as the beat design, and $K$ is the number of beats at the last iteration. For the predicted workload in month $\ell$, the greedy algorithm can be performed by selecting the beat $k$ in $D_n$ with the largest workload, i.e., $\arg\max_k \{w_{k\ell}(D_n)\}_{k \in \mathscr{K}}$. Then we split up the beat $k$ evenly into two beats using the K-means algorithm, where each atom in the beat is considered as a point. This will lead to generating a new beat, i.e., $K \coloneq K + 1$ and $\mathscr{K} \eqcolon \mathscr{K} \cup {K}$. The above process can be carried out iteratively until we find the design with the optimal number of beats.



We visualize our greedy design with different number of beats in Figures~\ref{fig:greedy-design}.
As seen from the result, the beat with the highest workload, shown in red, is split in each iteration as a result. We also examine the variance of beat workload versus number of beats, and find the optimal number of beats,
which will be further discussed in Sec.~\ref{sec:results}. 

\vspace{.05in}
\noindent{\bf Mixed-integer programming.}
Mathematical programming models 
are essential tools for modeling and solving redistricting problems, which can guarantee the optimality of the obtained solutions, are mostly based on mixed-integer programming (MIP). 
However, as shown in Sec.~\ref{sec:constraints}, the problem involves a large number of variables, including 21,134,535 continuous variables and 35,610 binary variables, as well as a set of additional linear constraints needed to be satisfied. 
In practice, the problem itself of searching for the global optimal design is computationally intractable and hard to be implemented on a large scale. 


\vspace{.05in}
\noindent{\bf Heuristic search.}
A metaheuristic method, simulated annealing (see, e.g., \cite{bertsimas1993simulated}), has been widely adopted in solving the large-scale combinatorial optimization problem. 
The simulated annealing algorithm explores the neighborhood of the current solution and decides a better substitution randomly. Simulated annealing can achieve reasonable performance in practice for various settings, although there are very limited theoretical performance guarantees \cite{aarts1987simulated, van1987simulated,lecchini2008simulated}. 
In particular, in our setting, we use the current/existing partition as an initial solution. Based on this, a new solution can be founded by selecting from a set of candidate solutions. The set of candidate solutions is typically constructed as ``neighboring'' solutions to the current solution without breaking contiguity. 

Specifically, in the $n$-th iteration, our simulated annealing algorithm performs the following acceptance-rejection sampling. Suppose the starting partition is $\mathcal{P}_n$. For instance, we can take the existing partition as an initialization. The next partition $\mathcal{P}_{n+1}$ is selected from a set of candidate partitions defined as $\mathcal S_{n+1}$ and $\mathcal{P}_{n+1} \in \mathcal{S}_{n+1}$.  
The candidate partitions in $\mathcal S_{n+1}$ satisfy contiguity and balance constraints. We randomly choose one of these candidate partitions $\mathcal{P}_{n+1} \in \mathcal S_{n+1}$, and evaluate a score   %
\[
    P(\mathcal{P}_{n+1}, \mathcal{P}_n|T) = 
    \begin{cases}
    1, & Z(\mathcal{P}_{n+1}) < Z(\mathcal{P}_n),\\
    \exp\{|Z(\mathcal{P}_{n
    +1}) - Z(\mathcal{P}_n)|/T\}, & \text{otherwise}.\\
    \end{cases}
\] 
where $Z(\cdot)$ denotes the cost associated with a partition (e.g., the compactness shown in \eqref{eq:opt-con-contiguity}), $T$ is a pre-specified temperature parameter that determines the speed of convergence, and the $P$ is the acceptance probability. We generate an independent uniform random variable $U \in [0, 1]$. 
The proposed partition is accepted if $P(\mathcal{P}_{n+1}, \mathcal{P}_n|T) \ge U$. We refer to an update of the proposed partition as a \emph{transition}. Note that there is a chance that the transition happens from a ``low-cost'' partition to a ``high-cost'' partition, and this ``perturbation'' will prevent the algorithm from being trapped at a local sub-optimal solution. The choice of the set candidate partitions $\mathcal S_{n+1}$ is critical for the performance of simulated annealing, which involves the trade-off between exploration and exploitation. Below, we introduce two strategies for candidate partitions and explore two types of ``neighbor'' partitions based on square and hexagonal grids, respectively.

\begin{figure}[!t]
    \includegraphics[width=1.\linewidth]{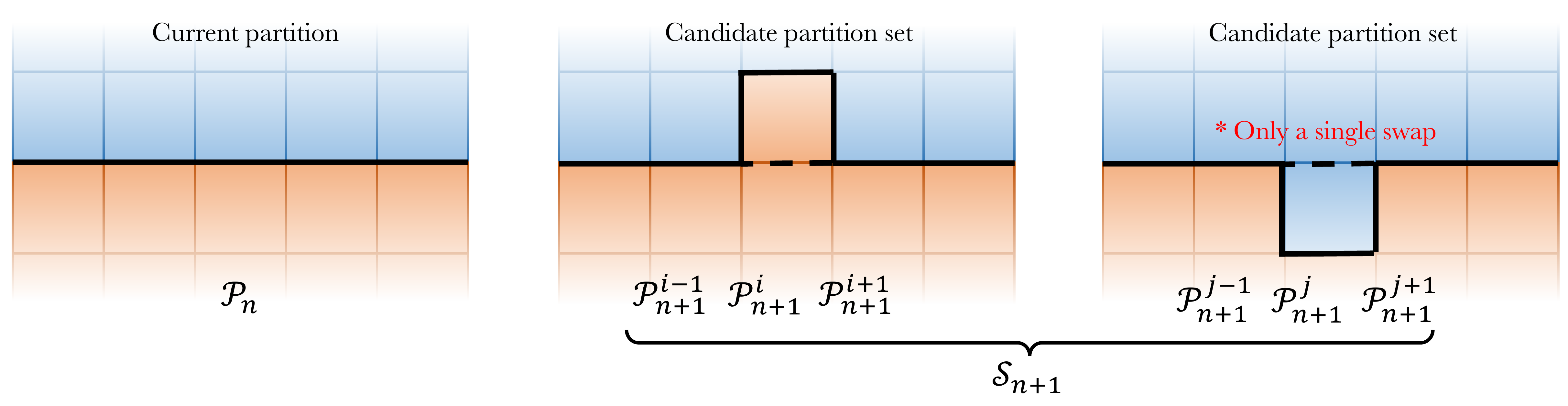}
    \vspace{-0.2in}
\caption{Illustrations of two approaches for candidate partitions based on one-swapping neighborhoods. Red and blue boxes represent vertices in different parts. The thick black line represents the boundary of two parts. The left panel shows the current partitions $\mathcal{P}_n$; the middle and the right panels show the candidate partition sets for the next iteration $\mathcal{S}_{n+1}$.}
\vspace{-0.1in}
\label{fig:sa-demo}
\end{figure}

As illustrated in Figure~\ref{fig:sa-demo}, we first consider the following simple heuristic in constructing the candidate set. This allows us to search for local optimal partitions at a reasonable computational cost~\cite{Wang2013}. The candidate set contains all partitions that swap a single vertex assignment at the boundary of the current partition. This simple heuristic is easy to implement since the number of such candidate partitions is usually small (because we only swap one end of the boundary edge). However, on the other hand, such candidate sets may contain partitions that are still too similar to the current partition. Therefore, we will consider the following alternative strategy.

\section{Results}
\label{sec:results}

In this section, we present our numerical results and final beat redesign for the City of South Fulton.

\begin{figure}[!t]
\centering
\begin{subfigure}[h]{0.32\linewidth}
\includegraphics[width=\linewidth]{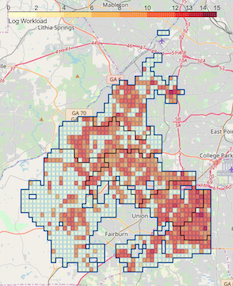}
\caption{2020 prediction}
\label{fig:2019 workload}
\end{subfigure}
\begin{subfigure}[h]{0.32\linewidth}
\includegraphics[width=\linewidth]{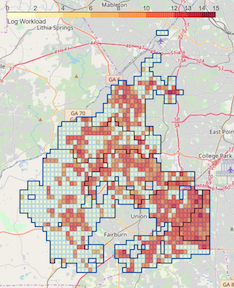}
\caption{2021 prediction}
\label{fig:2020 workload}
\end{subfigure}
\begin{subfigure}[h]{0.32\linewidth}
\includegraphics[width=\linewidth]{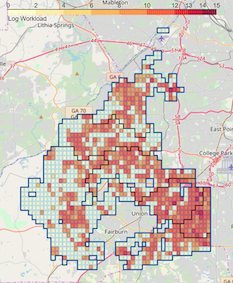}
\caption{2022 prediction}
\label{fig:2021 workload}
\end{subfigure}
\vspace{-0.1in}
\caption{Workload prediction where dark lines outline boundaries of beats and the color depth represents the level of the atom workload in each year.}
\label{fig:Workload prediction}
\vspace{-0.1in}
\end{figure}

\vspace{0.05in}
\noindent{\bf Workload analysis and prediction.}
The most important metric for evaluating imbalance we considered is workload variance over beats. As we defined in Sec.~\ref{sec:beat_redesign_opt}, the variance is the sum of the squared deviation of the beat workload from its mean. To fully understand the workload imbalance situation, it is necessary to show how the existing configuration exacerbates the unbalance of workload over beats in the past and how the existing configuration will impact the future. 
Figures~\ref{fig:Workload prediction} summarizes the predicted workload distribution over the entire city for the next four years from 2020 to 2022. As we can see from the map, there is a clear trend that the general workload level continues to increase, and the major workload concentrates on particular areas (such as College Park in the east of the city and I-285 \& I-20 in beat 4).
Due to the increasing growth of South Fulton and urban sprawl, this trend is leading to a police workload imbalance. 

\begin{figure*}
\centering
\begin{subfigure}[h]{0.16\linewidth}
\includegraphics[width=\linewidth]{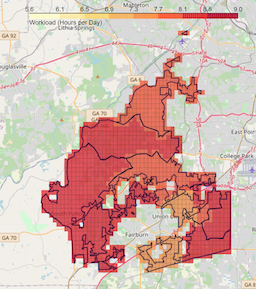}
\caption{existing 7-beat}
\label{fig:design-origin}
\end{subfigure}
\begin{subfigure}[h]{0.16\linewidth}
\includegraphics[width=\linewidth]{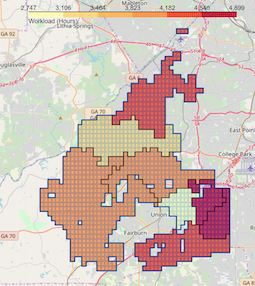}
\caption{greedy 8-beat}
\label{fig:8 beat design}
\end{subfigure}
\begin{subfigure}[h]{0.16\linewidth}
\includegraphics[width=\linewidth]{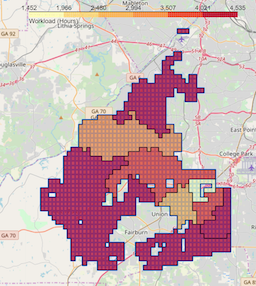}
\caption{greedy 9-design}
\label{fig:9 beat design}
\end{subfigure}
\begin{subfigure}[h]{0.16\linewidth}
\includegraphics[width=\linewidth]{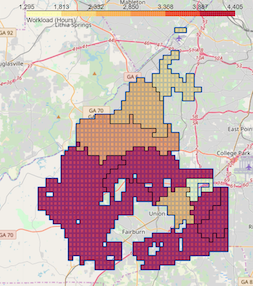}
\caption{greedy 10-beat}
\label{fig:10 beat design}
\end{subfigure}
\begin{subfigure}[h]{0.16\linewidth}
\includegraphics[width=\linewidth]{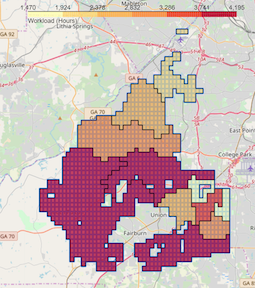}
\caption{greedy 11-beat}
\label{fig:11 beat design}
\end{subfigure}
\begin{subfigure}[h]{0.16\linewidth}
\includegraphics[width=\linewidth]{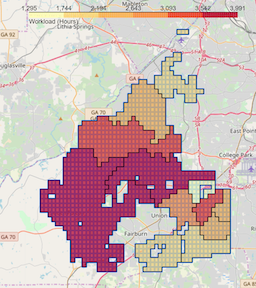}
\caption{greedy 12-beat}
\label{fig:12 beat design}
\end{subfigure}
\vfill
\begin{subfigure}[h]{0.16\linewidth}
\includegraphics[width=\linewidth]{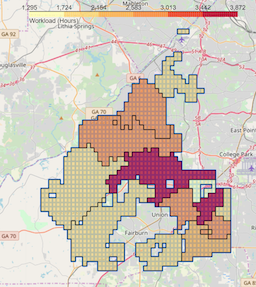}
\caption{greedy 13-beat}
\label{fig:13 beat design}
\end{subfigure}
\begin{subfigure}[h]{0.16\linewidth}
\includegraphics[width=\linewidth]{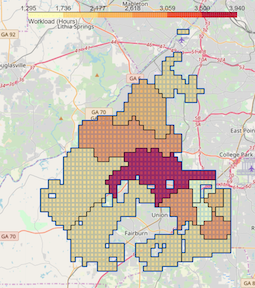}
\caption{greedy 14-beat}
\label{fig:14 beat design}
\end{subfigure}
\begin{subfigure}[h]{0.16\linewidth}
\includegraphics[width=\linewidth]{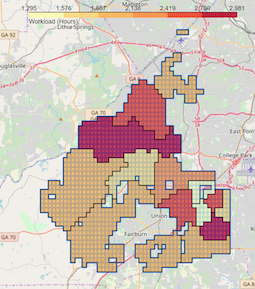}
\caption{greedy 15-beat}
\label{fig:15 beat design}
\end{subfigure}
\begin{subfigure}[h]{0.16\linewidth}
\includegraphics[width=\linewidth]{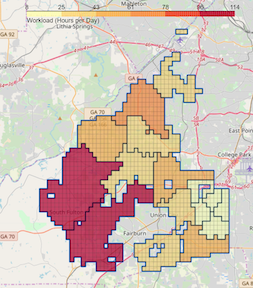}
\caption{greedy 16-beat}
\label{fig:16 beat design}
\end{subfigure}
\begin{subfigure}[h]{0.16\linewidth}
\includegraphics[width=\linewidth]{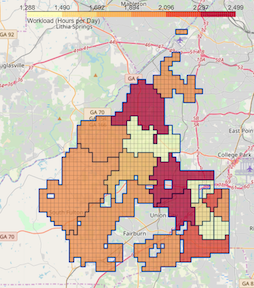}
\caption{greedy 17-beat}
\label{fig:17 beat design}
\end{subfigure}
\begin{subfigure}[h]{0.16\linewidth}
\includegraphics[width=\linewidth]{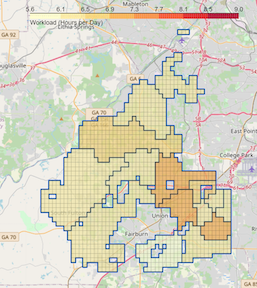}
\caption{proposed 15-beat}
\label{fig:15beat_design}
\end{subfigure}
\vspace{-0.1in}
\caption{Greedy beat designs where dark lines outline boundaries of beats and the color depth represents the level of the beat workload. The scale is adjusted in each image.}
\label{fig:greedy-design}
\vspace{-0.1in}
\end{figure*}

\begin{figure}[!b]
\vspace{-0.2in}
\centering
\includegraphics[width=.7\linewidth]{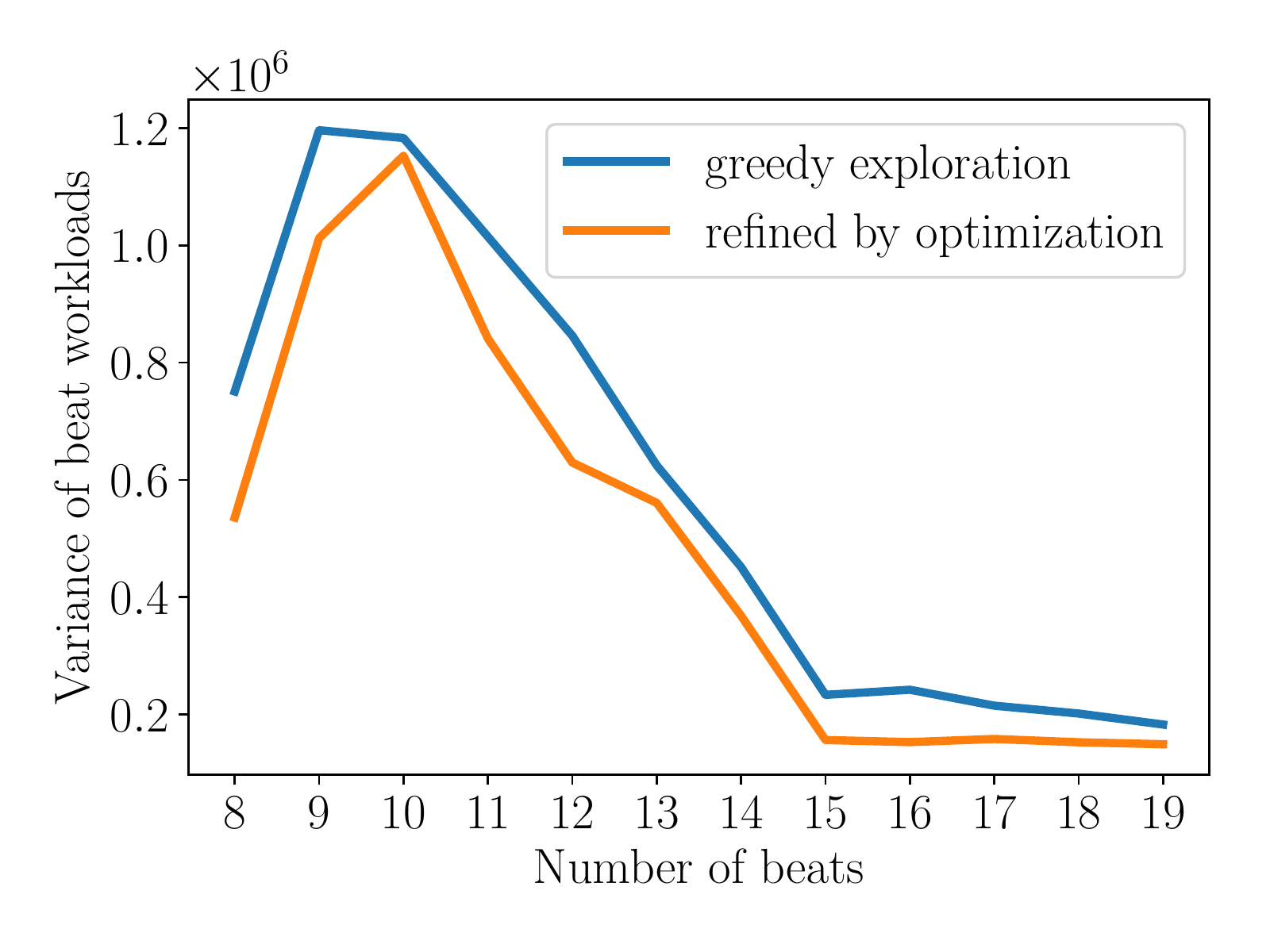}
\vspace{-0.1in}
\caption{Workload variance with different number of beats.}
\label{fig:optimal-num-beats}
\end{figure}

\vspace{0.05in}
\noindent{\bf Optimal beat number.}
When creating a beat design, the most important metric for evaluating imbalance is the workload variance over beats. However, for determining the optimal number of beats in the design, we also need to consider the cost associated with adding more beats, which includes the cost of additional training, hiring new officers, and so on. Therefore, there is a trade-off to minimize the workload variance while avoiding unnecessary costs for adding new beats.   
Figure~\ref{fig:greedy-design} presents comparisons between existing beat design, designs generated by greedy exploration algorithm, and the proposed design. 
Figure~\ref{fig:optimal-num-beats} shows that as we first begin to increase the number of beats, the workload variance decreases sharply before 15 beats.
We have shown that there are diminishing returns as we further increase the number of beats beyond 15. 
Therefore, we call 15 the optimal number of beats and the corresponding 15-beat greedy design will be used as an initialization of the simulated annealing for further refinement.




\begin{table}[!b]
\caption{Summarization of workload per beat.}
\label{tab:config-comparison}
\vspace{-0.15in}
\begin{center}
\begin{small}
\begin{sc}
\resizebox{0.4\textwidth}{!}{%
\begin{tabular}{cccccc}
\toprule[1pt]\midrule[0.3pt]
     Beat Number & \multicolumn{3}{c}{Workload in 2019} & \multicolumn{2}{c}{Workload in 2021}\\ 
       & \multicolumn{3}{c}{(hours/day)} & \multicolumn{2}{c}{(hours/day)}\\ 

     & Existing & Greedy & Refined & Greedy & Refined \\
     \hline
     1 & 38.59 & 17.15 & 17.15 & 18.05 & 18.05\\
     2 & 24.84 & 24.84 & 23.56 & 27.09 & 25.61\\
     3 & 32.84 & 18.78 & 20.08 & 17.91 & 19.91\\
     4 & 34.44 & 17.45 & 17.08 & 16.83 & 16.14 \\
     5 & 65.94 & 22.10 & 20.31 & 21.40 & 19.32 \\
     6 & 38.44 & 14.69 & 18.30 & 14.54 & 16.73 \\ 
     7 & 34.96 & 17.55 & 19.99 & 17.67 & 20.01 \\
     8 & N/A & 12.51 & 12.51 & 11.66 & 11.66 \\
     9 & N/A & 10.79 & 10.79 & 11.10 & 11.10 \\
     10 & N/A & 21.45 & 21.87 & 21.45 & 21.87 \\
     11 & N/A & 23.75 & 19.33 & 22.2 & 22.62 \\
     12 & N/A & 17.41 & 17.41 & 23.40 & 21.60 \\
     13 & N/A & 17.00 & 16.82 & 16.70 & 15.87 \\
     14 & N/A & 20.53 & 18.89 & 19.99 & 17.81 \\
     15 & N/A & 14.06 & 15.94 & 13.18 & 15.93 \\
     Variance & 142.91 & 15.12 & 10.13 & 18.269 & 13.15 \\
\midrule[0.3pt]\bottomrule[1pt]
\end{tabular}
}
\end{sc}
\end{small}
\end{center}
\vspace{-0.1in}
\end{table}

\vspace{0.05in}
\noindent{\bf Proposed beat design.}
The initial report in 2019 contained beat-wise workload prediction for the next three years (2020, 2021 and 2022), and proposed three candidate designs with similar beat shifts that all attains the best workload balance. 
In Table~\ref{tab:config-comparison}, we list the predicted annual workload in each beat, total workload, and workload variance. 
After the plans were reported to the police, we met several times to deliberate the various trade-offs, and held police-engagement meetings to elicit feedback from the patrol force. Police Deputy Chief and couples of key senior officers also participated in these discussions and voiced comments.
The new design and the previous existing design have been both presented in Figure~\ref{fig:greedy-design}.
This design is preferred by the police for three major reasons: (2) this plan makes the minimal changes based on the existing police zone configuration in comparison with other candidate plans, which minimizes the implementation cost in practice; (2) the workload variance has been drastically reduced by 89\% $\sim$ 92\% by increasing the number of beats to 15; (3) the proposed plan achieves a lower level of workload variance as well as a smaller variance increment in the future year 2021.

\vspace{0.05in}
\noindent{\bf Staffing level analysis.}
We quantify our potential police response workload by converting the workload in each beat into hours per day. Table~\ref{tab:config-comparison} shows real workload distribution in 2019 and predicted workload distribution in 2021 under different designs, respectively. Entries of the table suggest the number of hours per day, a police officer would expect to be responding to 911 calls. As we can see, our proposed beat design drastically reduces the beatwise workload. In particular, the proposed design results in a decrease in workload variance of over 85\% comparing to the existing design, making policing more equitable in the city.

In the City of South Fulton council meeting, the city council emphasizes the importance of community engagement from the police force. Thanks to our beat design, the police workload per day in each beat can be reduced drastically; this will allow police officers to participate in community events and start pro-active patrols. This is a huge difference from the past 50 years, where police officers have been going from call to call on their entire shift. Additionally, the staffing level prediction gives the SFPD how many officers they need to handle the 911 calls in a beat. They then can recruit more officers for the sole purpose of community engagement and pro-active patrolling if they desire. 

\section{Implementation}

In January 2020, we submitted the final report to the South Fulton Police Department and the South Fulton City Council. The report was reviewed by Police Chief Meadows, Deputy Police Chief Rogers, and Mayor Bill Edwards. Our report analyzed the police workload and proposed a detailed redistricting plan. Our redistricting plan mainly changed in four areas (Figure~\ref{fig:15beat_design}): We add three new beats in the southeast of the city near College Park, the area with the highest workload. The biggest beat in the west of the city is split into two beats. We add a beat in the north of the city near the airport. The southern beat is also split into two. In total, the redistricting plan has reduced the response time throughout the city and rebalanced the police workload between the fifteen beats.

Later that month, the South Fulton City Council approved the new beat design. The South Fulton Police Department plan to implement the new beat design in early 2020. 
The new beat design was praised by the city council, as some council members said that our beat design and study has been long needed and that it sets an example for other cities in the southeast. 
Residents of South Fulton acclaimed about the change on social media and thanked the City of South Fulton Police Department and our team for contributing to the communities. The new beat design also received coverage from several news sources, including Fox 5 Atlanta \cite{FoxNews}.



\section{Conclusion}

In this paper, we presented our work on the City of South Fulton police beat redesign. We propose an optimization framework with the spatial regression model as well as large-scale data analytics. We construct an operational model to predict zone workload using an accurate and tractable linear approximation. The proposed method yields a redesigned zone plan with lower workload variance by only changing eight beats.  Currently, we are continuing our partnership with the SFPD. We will continue to observe the police workload in the City of South Fulton as the city and workload grow. If the workload becomes unbalanced once more, we can quickly suggest a new beat design using our already existing methods. As the SFPD continues to grow, they will also hire an information officer that will assist in workload analytics and carry on our workload prediction.



\newpage

\section*{Acknowledgments}
 This work is supported in part by the National Science Foundation under Grant CMMI-2015787
\bibliographystyle{ACM-Reference-Format}
\bibliography{refs}










\end{document}